\documentclass[11pt]{article}

\newcommand{\real}{{\bf R}}

\newcommand{\intplus}{{\bf N}}

\renewcommand{\div}{\mathop{\rm div}}

\newcommand{\e}{{\rm e}}                
\renewcommand{\d}{\,{\rm d}}            
\newcommand{\D}{{\rm d}}                

\newcommand{\meas}{{\rm meas}}

\newcommand{\cH}{{\cal H}}

\newcommand{\cM}{{\cal M}}

\newtheorem{theorem}{Theorem}[section]

\newtheorem{definition}[theorem]{Definition}
\newtheorem{proposition}[theorem]{Proposition}

\newtheorem{remark}[theorem]{Remark}

\newcommand{\reff}[1]{(\ref{#1})}

\newcommand{\inttwo}{\int_{\real^2}}
\newcommand{\proof}{{\noindent \bf Proof:\ }}

\def\epsilon{\varepsilon}
\def\phi{\varphi}
\def\weakto{\rightharpoonup}

\def\build#1_#2^#3{\mathrel{
  \mathop{\kern 0pt#1}\limits_{#2}^{#3}}}
\def\QED{\mbox{}\hfill$\Box$}

\newcommand{\one}{\mathbf{1}}

\usepackage{amsfonts}
\begin{document}

\title{On the uniqueness of the solution of the two-dimensional 
  Navier-Stokes equation with a Dirac mass as initial vorticity}

\author{
Isabelle Gallagher \\Universit{\'e} Paris 7
\\
Institut de Math{\'e}matiques de Jussieu \\
Case 7012, 2 place Jussieu
 \\
75251 Paris Cedex 05
 France 
\and
Thierry Gallay \\ 
Institut Fourier \\
Universit{\'e} de Grenoble I \\
38402 Saint-Martin d'H{\`e}res, France 
\vspace{1mm}
\and
Pierre-Louis Lions \\
CEREMADE \\
Universit{\'e} de Paris-Dauphine \\
75775 Paris cedex 16, France}
\date{}
\maketitle
\begin{abstract}
We propose two different proofs of the fact that Oseen's vortex is 
the unique solution of the two-dimensional Navier-Stokes equation
with a Dirac mass as initial vorticity. The first argument, due to
C.E.~Wayne and the second author, is based on an entropy estimate
for the vorticity equation in self-similar variables. The second
proof is new and relies on symmetrization techniques for 
parabolic equations.
\end{abstract}


\section{Introduction}\label{intro} 
\setcounter{equation}{0}

We consider the vorticity equation associated to the two-dimensional
Navier-Stokes equation, namely
\begin{equation}\label{Veq}
  \partial_t \omega(x,t) + 
  u(x,t) \cdot \nabla \omega(x,t)  \,=\, \Delta \omega(x,t)~, 
  \quad x \in \real^2~, \quad t > 0~.
\end{equation}
The velocity field $u(x,t) \in \real^2$ is obtained from the vorticity
$\omega(x,t) \in \real$ via the Biot-Savart law
\begin{equation}\label{BS}
  u(x,t) \,=\, \frac{1}{2\pi} \inttwo 
  \frac{(x - y)^{\perp}}{|x-y|^2}\, \omega(y,t)\d y\ ,
  \quad x \in \real^2~, \quad t > 0~,
\end{equation}
where $(x_1,x_2)^\perp = (-x_2,x_1)$. It satisfies $\div u = 0$
and $\partial_1 u_2 - \partial_2 u_1 = \omega$. Equations 
\reff{Veq}, \reff{BS} are invariant under the scaling transformation
\begin{equation}\label{scaling}
  \omega(x,t) \,\mapsto\, \lambda^2 \omega(\lambda x,\lambda^2 t)~, \quad 
  u(x,t) \,\mapsto\, \lambda u(\lambda x,\lambda^2 t)~, \quad 
  \lambda > 0~.
\end{equation}

The Cauchy problem for the vorticity equation \reff{Veq} is globally
well-posed in the (scale invariant) Lebesgue space $L^1(\real^2)$, see
for instance \cite{ben-artzi}. To include more general initial data,
such as isolated vortices or vortex filaments, it is necessary to use
larger function spaces. A natural candidate is the space
$\cM(\real^2)$ of all finite real measures on $\real^2$, equipped with
the total variation norm. This space contains $L^1(\real^2)$ as a
closed subspace, and its norm is invariant under (the spatial part of)
the rescaling \reff{scaling}. Moreover, $\cM(\real^2)$ is closed with
respect to the following weak convergence: $\mu_n \weakto \mu$ if
$\inttwo \phi \d\mu_n \to \inttwo \phi \d\mu$ for any continuous
function $\phi : \real^2 \to \real$ vanishing at infinity.

Existence of solutions of \reff{Veq} with initial data in 
$\cM(\real^2)$ was first proved by Cottet \cite{cottet}, 
and independently by Giga, Miyakawa and Osada \cite{GMO}, 
see also Kato \cite{kato}. Uniqueness can be obtained by a 
standard Gronwall argument if the {\em atomic part} of the 
initial vorticity $\mu$ is sufficiently small \cite{GMO,kato}, 
but this method is bound to fail if $\mu$ contains large 
Dirac masses. In the particular case where $\mu = \alpha \delta_0$
for some $\alpha \in \real$, an explicit solution is known:
\begin{equation}\label{oseendef}
  \omega(x,t) \,=\, \frac{\alpha}{t}\,G\Bigl(
    \frac{x}{\sqrt{t}}\Bigr)~, \quad
  u(x,t) \,=\, \frac{\alpha}{\sqrt{t}}\,v^G\Bigl(
    \frac{x}{\sqrt{t}}\Bigr)~, \quad x \in \real^2~, \quad t > 0~,
\end{equation}
where
\begin{equation}\label{GvGdef}
   G(\xi) \,=\, \frac{1}{4\pi}\,\e^{-|\xi|^2/4}~, \quad
   v^G(\xi) \,=\, \frac{1}{2\pi}\,\frac{\xi^\perp}{|\xi|^2}
   \,\Bigl(1 - \e^{-|\xi|^2/4}\Bigr)~, \quad \xi \in \real^2~.
\end{equation}
This self-similar solution of the two-dimensional Navier-Stokes
equation is often called the {\em Lamb-Oseen vortex} with total
circulation $\alpha$. It is the unique solution with initial 
vorticity $\alpha \delta_0$ in the following precise sense:

\begin{theorem}\label{main} {\bf \cite{gallay-wayne}}
Let $T > 0$, $K > 0$, $\alpha \in \real$, and assume that
$\omega \in C^0((0,T),L^1(\real^2)\cap L^\infty(\real^2))$
is a solution of \reff{Veq} satisfying $\|\omega(\cdot,t)\|_{L^1}
\le K$ for all $t \in (0,T)$ and $\omega(\cdot,t) \weakto \alpha
\delta_0$ as $t \to 0+$. Then
$$
   \omega(x,t) = \frac{\alpha}{t} \,G\Bigl(\frac{x}{\sqrt{t}}
   \Bigr)~, \quad x \in \real^2~,\quad t \in (0,T)~.
$$
\end{theorem}

Here and in the sequel, we say that $\omega \in C^0((0,T),L^1(\real^2)
\cap L^\infty(\real^2))$ is a (mild) solution of~\reff{Veq} if the 
associated integral equation
\begin{equation}\label{Vint}
  \omega(\cdot,t_2) \,=\, \e^{(t_2-t_1)\Delta}\omega(\cdot,t_1) 
   - \int_{t_1}^{t_2} \nabla \cdot \e^{(t_2-t)\Delta} 
  \Bigl(u(\cdot,t)\omega(\cdot,t)\Bigr)\d t
\end{equation}
is satisfied for all $0 < t_1 < t_2 < T$.

The particular case covered by Theorem~\ref{main} is important for at
least two reasons. First, it has deep connections with the long-time
behavior of smooth solutions. Indeed, if $\omega(x,t)$ is any solution
of \reff{Veq} such that $\omega(\cdot,t) \in L^1(\real^2)$ and
$\inttwo \omega(x,t)\d x = \alpha$ for all $t > 0$, it is shown in
\cite{gallay-wayne} that~$\omega(\cdot,t)$ converges in $L^1(\real^2)$ to
the Oseen vortex with circulation $\alpha$ as $t \to +\infty$.  This
convergence result can in fact be {\em deduced} from
Theorem~\ref{main} by a classical renormalization argument, see
\cite{carpio} and~\cite{gigabook}. On the other hand, combining the
standard Gronwall argument (which works if the initial measure has
small atomic part) with Theorem~\ref{main} (which allows to treat a
large Dirac mass), it is possible to prove the uniqueness of the
solution of \reff{Veq} for {\em arbitrary data} $\mu \in \cM(\real^2)$
\cite{GG}. Thus, the Cauchy problem for \reff{Veq} is globally
well-posed in the space $\cM(\real^2)$ too.

The proof of Theorem~\ref{main} in \cite{gallay-wayne} is based
on an entropy estimate for equation \reff{Veq} rewritten in 
self-similar variables. For completeness, this argument will 
be reproduced below in its simplest form. The real purpose of
this paper is to provide an alternative proof of Theorem~\ref{main}
which relies on completely different ideas: it is based on 
symmetrization techniques for elliptic and parabolic equations 
which originate in the work of Talenti \cite{talenti}, see 
\cite{ALT1,ALT2}. These methods ultimately rely on the maximum 
principle, whereas the original proof of \cite{gallay-wayne}
uses a Lyapunov function. We find the comparison of both arguments
very instructive, and this is why we chose to collect them in 
the present article. 

The rest of the text is organized as follows. Under the assumptions of
Theorem~\ref{main}, we first derive a priori estimates on the solution
$\omega(x,t)$ using Gaussian bounds on the fundamental solution of the
convection-diffusion equation, which are due to Osada in~\cite{osada}. 
This preliminary part is common to both approaches. Following 
\cite{gallay-wayne}, we next introduce the self-similar variables~$\xi = x/\sqrt{t}$,~$\tau = \log(t)$ and we show that the rescaled 
vorticity $w(\xi,\tau) = \e^\tau \omega(\xi\,\e^{\frac\tau2},\e^\tau)$ 
satisfies an entropy estimate which implies the conclusion of
Theorem~\ref{main}. The new proof begins in
Section~\ref{symmetrization}, where we recall the definition of the
symmetric nonincreasing rearrangement of a function. In particular, we
prove that the solution of a convection-diffusion equation of the
form~\reff{Veq} 
with initial data $\omega_0$ is ``dominated'' (in a sense
to be precised) for all times by the solution of the {\em heat
equation} with symmetrized initial data $\omega_0^\#$. In the final
section, we show how this key estimate provides a new and short proof
of Theorem~\ref{main}.


\section{Representation and a priori estimates}
\label{representation}
\setcounter{equation}{0}

Assume that $\omega \in C^0((0,T),L^1(\real^2) \cap 
L^\infty(\real^2))$ is a solution of \reff{Veq} satisfying 
the hypotheses of Theorem~\ref{main}. From~\cite{brezis} we know 
that~$\omega(x,t)$ coincides for $t > 0$ with a classical solution 
of \reff{Veq} in~$\real^2$ as constructed for instance
in~\cite{ben-artzi}. In particular~$\omega(x,t)$ is 
smooth for~$t > 0$.
Since the Cauchy problem for \reff{Veq} is globally well-posed
in $L^1(\real^2) \cap L^\infty(\real^2)$ and since the~$L^1$ norm
of $\omega(\cdot,t)$ is nonincreasing with time, we can assume that~$T = +\infty$ without 
loss of generality. By assumption $\|\omega
(\cdot,t)\|_{L^1} \le K$ for all $t > 0$, hence the associated 
velocity field satisfies~$t^{\frac12}\|u(\cdot,t)\|_{L^\infty} 
\le CK$ for all~$t > 0$, where~$C > 0$ is a universal constant, 
see (\cite{carlen-loss}, Theorem~2).

The solution $\omega(x,t)$ of \reff{Veq} has an integral representation 
of the form
\begin{equation}\label{omrepr}
  \omega(x,t) \,=\, \inttwo \Gamma_u(x,t;y,s)\omega(y,s)\d y~,
  \quad x \in \real^2~, \quad t > s > 0~,
\end{equation}
where $\Gamma_u$ is the fundamental solution of the 
convection-diffusion equation $\partial_t \omega + u \cdot
\nabla \omega = \Delta \omega$. The following properties of $\Gamma_u$
are due to Osada \cite{osada} and to Carlen and Loss
\cite{carlen-loss}.

\medskip\noindent{$\bullet$} For any $\beta \in (0,1)$
there exists $K_1 > 0$ (depending only on $K$ and $\beta$) 
such that 
\begin{equation}\label{gamm1}
  0 \,<\, \Gamma_u(x,t;y,s) \,\le\, \frac{K_1}{t-s}
  \,\exp\Bigl(-\beta\frac{|x-y|^2}{4(t-s)}\Bigr)\ ,
\end{equation}
for all $x,y \in \real^2$ and $t > s > 0$, see \cite{carlen-loss}. 
A similar Gaussian lower bound is also known. 

\medskip\noindent{$\bullet$} There exists $\gamma \in (0,1)$ 
(depending only on $K$) and, for any $\delta > 0$, there exists
$K_2 > 0$ (depending only on $K$ and $\delta$) such that 
\begin{equation}\label{gamm2}
  |\Gamma_u(x,t;y,s) - \Gamma_u(x',t';y',s')| \le K_2
  \Bigl(|x{-}x'|^\gamma + |t{-}t'|^{\gamma/2} + 
  |y{-}y'|^\gamma + |s{-}s'|^{\gamma/2}\Bigr)\ ,
\end{equation}
whenever $t-s \ge \delta$ and $t'-s' \ge \delta$, see \cite{osada}.

\medskip\noindent{$\bullet$} For $t > s > 0$ and $x,y \in
\real^2$,  
\begin{equation}\label{gamm3}
  \inttwo \Gamma_u(x,t;y,s)\d x \,=\, 1~, \quad 
  \inttwo \Gamma_u(x,t;y,s)\d y \,=\, 1~. 
\end{equation}

\medskip\noindent 
If $x,y \in \real^2$ and $t > 0$, it follows from \reff{gamm2} that 
the function $s \mapsto \Gamma_u(x,t;y,s)$ can be continuously
extended up to $s = 0$, and that this extension (still denoted by 
$\Gamma_u$) satisfies properties~\reff{gamm1} to \reff{gamm3} with 
$s = 0$. Using this observation, we obtain for all $x \in \real^2$ and 
all $t > 0$ 
\begin{eqnarray*}
  \omega(x,t) &=& \inttwo \Gamma_u(x,t;y,0)\omega(y,s)\d y \\
  &+& \inttwo \Bigl(\Gamma_u(x,t;y,s) - \Gamma_u(x,t;y,0)\Bigr)
  \omega(y,s)\d y~, \quad 0 < s < t~.
\end{eqnarray*}
Since $\omega(\cdot,s)$ is bounded in $L^1(\real^2)$, it follows
from \reff{gamm2} that the second integral in the right-hand side
converges to zero as $s$ goes to zero. On the other hand, since $y
\mapsto \Gamma_u(x,t;y,0)$ is continuous and vanishes at infinity, and
since $\omega(\cdot,s) \weakto \alpha\delta_0$ as $s \to 0$, we can 
take the limit $s \to 0$ in the first integral and we obtain the 
useful representation:
\begin{equation}\label{repres}
  \omega(x,t) \,=\, \alpha \Gamma_u(x,t;0,0)~, \quad  
  \quad x \in \real^2~, \quad t > 0~.
\end{equation}

This formula shows in particular that $\omega \equiv 0$ if 
$\alpha = 0$. Thus, upon replacing~$\omega(x_1,x_2,t)$ with~$-\omega(x_2,x_1,t)$ if needed, we can assume from now on that 
$\alpha > 0$. To simplify the notation, we denote by $\Omega(x,t)$ 
the Oseen vortex with circulation $\alpha$, namely
\begin{equation}\label{Omegadef}
  \Omega(x,t) \,=\, \frac{\alpha}{4\pi t}\,\e^{-\frac{|x|^2}{4t}}~,
  \quad x \in \real^2~, \quad t > 0~.
\end{equation}
From \reff{gamm1}, \reff{repres} we have
\begin{equation}\label{gaussbound}
  0 \,<\, \omega(x,t) \,\le\, \frac{K_1 \alpha}{t}
  \,\e^{-\beta\frac{|x|^2}{4t}}~, \quad x \in \real^2~, \quad t > 0~.
\end{equation}
Moreover, in view of \reff{gamm3}, 
\begin{equation}\label{mom0}
  \inttwo \omega(x,t)\d x \,=\, \alpha \,=\, \inttwo \Omega(x,t)
  \d x~, \quad t > 0~.
\end{equation}
In particular $\|\omega(\cdot,t)\|_{L^1} = \alpha$ for all $t > 0$, 
so that we can take $K = |\alpha|$ in the statement of 
Theorem~\ref{main} without loss of generality. Finally, 
using \reff{Veq} and integrating by parts, we find
$$
   \frac{\D}{\D t} \inttwo |x|^2 \omega \d x \,=\, 
   \inttwo |x|^2 (\Delta \omega - u \cdot\nabla \omega)\d x 
   \,=\, 4\alpha + 2 \inttwo (x\cdot u)\omega \d x~.
$$
Actually, the last integral vanishes identically. This can 
be seen by replacing $u$ with its expression~\reff{BS} and using 
the symmetry properties of the 
Biot-Savart kernel. Thus we find that~$\frac{\D}{\D t} \inttwo |x|^2 \omega \d x 
= 4\alpha$, and in view of \reff{gaussbound} we conclude that
\begin{equation}\label{mom2}
  \inttwo |x|^2 \omega(x,t)\d x \,=\, 4 \alpha t \,=\, 
  \inttwo |x|^2 \Omega(x,t) \d x~, \quad t > 0~.
\end{equation}
 

\section{Uniqueness proof using entropy estimates}
\label{unique1}
\setcounter{equation}{0}

We first introduce the entropy functionals that will be used in the
proof. Let $f : \real^2 \to (0,+\infty)$ be a $C^1$ function 
satisfying $\inttwo f(\xi)\d \xi = 1$. We define
\begin{equation}\label{HIdef} 
  H(f) \,=\, \inttwo f(\xi) \log \Bigl(\frac{f(\xi)}{G(\xi)}\Bigr) 
  \d\xi~, \quad
  I(f) \,=\, \inttwo f(\xi) \Big|\nabla\log\Bigl(\frac{f(\xi)}{G(\xi)} 
  \Bigr)\Big|^2 \d\xi~,
\end{equation}
where $G$ is as in \reff{GvGdef}. In kinetic theory, $H(f)$ is 
the {\em relative Boltzmann entropy} of the distribution function 
$f$ with respect to the Gaussian $G$. In information theory, 
$H(f)$ is often called the {\em relative Kullback entropy} of $f$ 
with respect to $G$, and $I(f)$ the {\em relative Fisher information} 
of $f$ with respect to $G$. The entropy $H(f)$ satisfies the 
following bounds (see for instance \cite{AMTU})
\begin{equation}\label{Hineq}
  \frac12 \|f-G\|_{L^1}^2 \,\le\, H(f) \,\le\, I(f)~,
\end{equation}
where the first inequality is the Csisz{\'a}r-Kullback inequality
and the second bound is the Stam-Gross logarithmic Sobolev 
inequality. It follows in particular from \reff{Hineq} that
$H(f) \ge 0$, and~$H(f) = 0$ if and only if~$f = G$. 

Now, let $\omega(x,t)$ be a solution of \reff{Veq} satisfying 
the hypotheses of Theorem~\ref{main}, and let~$u(x,t)$ be the 
associated velocity field. As in the previous section, we assume 
without loss of generality that $T = +\infty$ and that $\alpha = 
\inttwo \omega(x,t)\d x > 0$. Following \cite{gallay-wayne}, we 
define the rescaled vorticity~$w(\xi,\tau)$ and the rescaled 
velocity field $v(\xi,\tau)$ by
\begin{equation}\label{wvdef}
  w(\xi,\tau) \,=\, \e^\tau \omega(\xi\,\e^{\frac\tau2},\e^\tau)~,
  \quad   
  v(\xi,\tau) \,=\, \e^{\frac\tau2} u(\xi\,\e^{\frac\tau2},\e^\tau)~,
  \quad \xi \in \real^2~, \quad \tau \in \real~.
\end{equation}
The equation satisfied by $w(\xi,\tau)$ reads
\begin{equation}\label{SV}
  \partial_\tau  w + (v \cdot \nabla_\xi) w \,=\, \Delta_\xi w 
  + \frac12 (\xi \cdot \nabla_\xi)w + w~.
\end{equation}
Since $\div v = 0$ and $\partial_1 v_2 - \partial_2 v_1 = w$, 
the rescaled velocity $v$ is obtained from the rescaled vorticity 
$w$ via the usual Biot-Savart law. 

Our aim is to show that $w(\xi,\tau) = \alpha G(\xi)$ for all 
$\xi \in \real^2$ and all $\tau \in \real$. To do that, we consider
the relative entropy $h(\tau) = H(w(\cdot,\tau)/\alpha)$, where 
$H$ is defined in \reff{HIdef}. Since by~\reff{gaussbound} we have~$0 < w(\xi,\tau) \le 
K_1 \alpha \,\e^{-\beta|\xi|^2/4}$, there
exists $K_3 > 0$ such that $0 \le h(\tau) \le K_3$ for all 
$\tau \in \real$. On the other hand, using \reff{SV} and integrating
by parts, it is not difficult to verify that
$$
   \frac{\D}{\D\tau} H\Bigl(\frac{w(\cdot,\tau)}{\alpha}\Bigr)
   \,=\, - I\Bigl(\frac{w(\cdot,\tau)}{\alpha}\Bigr)~, \quad 
   \tau \in \real~,
$$
see \cite{gallay-wayne}. In view of the second inequality in 
\reff{Hineq}, this implies that $h'(\tau) \le -h(\tau)$ for all 
$\tau \in \real$. Thus
$$
   h(\tau) \,\le\, \e^{-(\tau-\tau_0)} h(\tau_0) \,\le\, 
   K_3 \,\e^{-(\tau-\tau_0)} \quad \hbox{whenever}~ \tau \ge \tau_0~.
$$
Letting $\tau_0 \to -\infty$, we obtain $h(\tau) = 0$ for all 
$\tau \in \real$, hence $w(\xi,\tau) = \alpha G(\xi)$ for all 
$\xi \in \real^2$ and all $\tau \in \real$. This is the desired 
result. \QED


\section{Symmetric rearrangements and applications}
\label{symmetrization}
\setcounter{equation}{0}

In this section, we recall the definition of the symmetric 
nonincreasing rearrangement of a function $f : \real^N \to 
\real$, and we establish two results that will be used in the
final section. We refer to \cite{ALT0}, \cite{bandle}, 
\cite{lieb-loss} for more information on rearrangements and
their applications. 

Let $f : \real^N \to \real$ be a measurable function. We assume that 
$f$ vanishes at infinity in the sense that $\meas(\{x \in \real^N 
\,|\, |f(x)| > t\})  < \infty$ for all $t > 0$. We define the
distribution function~$\mu_f : [0,+\infty) \to [0,+\infty]$ by
$$
   \mu_f(t) \,=\, \meas(\{x \in \real^N \,|\, |f(x)| > t\})~, \quad
   t \ge 0~,
$$
and the nonincreasing rearrangement $f^* : [0,+\infty) \to [0,+\infty]$
by
$$
   f^*(s) \,=\, \sup\{t \ge 0 \,|\, \mu_f(t) > s\}~, \quad
   s \ge 0~.
$$
The {\em symmetric nonincreasing rearrangement} $f^\# : 
\real^N \to [0,+\infty]$ is then defined by
$$
   f^\#(x) \,=\, f^*(c_N |x|^N)~, \quad x \in \real^N~,
$$
where $c_N = \pi^{N/2}/\Gamma(\frac{N}2+1)$ is the measure of the
unit ball in $\real^N$. By construction, $f^\#$ is radially
symmetric and nonincreasing (along rays). Moreover, $f^\#$ is
lower semicontinuous (hence measurable) and $\meas(\{x \in \real^N 
\,|\, f^\#(x) > t\}) = \meas(\{x \in \real^N \,|\, |f(x)| > t\})$ 
for all $t > 0$. As a consequence, if $f \in L^p(\real^N)$ for
some $p \in [1,\infty]$, then $f^\# \in L^p(\real^N)$ and 
$\|f^\#\|_{L^p} = \|f\|_{L^p}$. More generally, one has 
$\|f^\#-g^\#\|_{L^p} \le \|f-g\|_{L^p}$ for all $f,g \in
L^p(\real^N)$. Note also that $f^\# = |f|^\#$, and that 
$f^\#$ is a continuous function if $f$ is continuous.  

We next introduce a partial order on integrable functions 
which is based on rearrangements.  

\begin{definition}\label{defdom}
Let $f,g : \real^N \to \real$ be integrable functions. 
We say that $f$ is {\em dominated} by $g$ if
\begin{equation}\label{domdef}
  \int_{B_R} f^\#(x)\d x \,\le\, \int_{B_R} g^\#(x)\d x
  \quad \hbox{for all } R > 0~,
\end{equation}
where $B_R = \{x \in \real^N \,|\, |x| < R\}$. In this case, 
we write $f \preceq g$.
\end{definition}

It can be shown that $f \preceq g$ if and only if
$$
   \int_{\real^N} \Phi(|f(x)|)\d x \,\le\, 
   \int_{\real^N} \Phi(|g(x)|)\d x
$$
for all convex functions $\Phi : [0,+\infty) \to [0,+\infty)$ with 
$\Phi(0) = 0$, see \cite{ALT0}. In particular, if $f \preceq g$
and $g \in L^1(\real^N) \cap L^\infty(\real^N)$, then $\|f\|_{L^p} 
\le \|g\|_{L^p}$ for any $p \ge 1$. 

The uniqueness proof in Section~\ref{unique2} will be based on 
two properties of the domination relation which we now describe. 
These results are rather standard, and the proofs will be outlined
for completeness only. 

\begin{proposition}\label{prop1}
Let $f,g : \real^N \to [0,+\infty)$ be continuous and integrable 
functions satisfying:\\[1mm]
\null\quad a) $f \preceq g$; \\[1mm]
\null\quad b) $g = g^\#$;\\[1mm]
\null\quad c) $\int_{\real^N} f(x)\d x = \int_{\real^N} g(x)\d x$;\\[1mm]
\null\quad d) $\int_{\real^N} |x|^N f(x)\d x = \int_{\real^N} |x|^N 
g(x)\d x < \infty$.\\[1mm]
Then $f = g$. 
\end{proposition}

\begin{remark}\label{moregeneral}
As is clear from the proof, the same result holds if $|x|^N$
is replaced by $\phi(|x|)$ in assumption d), where $\phi$ 
is any continuous, positive, and strictly increasing function
on $[0,+\infty)$. 
\end{remark}

\proof
We first observe that
\begin{equation}\label{concin}
  \int_{\real^N} |x|^N f^\#(x)\d x \,\le\, \int_{\real^N} 
  |x|^N f(x)\d x~,
\end{equation}
and that equality holds in \reff{concin} if and only if 
$f = f^\#$. Indeed, this is obvious if $f$ is the characteristic 
function of a bounded open set $A \subset \real^N$, because 
then $f^\# = \one_{B_R}$ where $c_N R^N = \meas(A)$. The general
case follows using the ``layer cake representation'' 
\cite{lieb-loss}
$$
   f(x) \,=\, \int_0^\infty \one_{\{f > t\}}(x) \d t~, \quad
   f^\#(x) \,=\, \int_0^\infty \one_{\{f > t\}}^\#(x) \d t~,
$$
together with Fubini's theorem. 

Now, let $h = g^*-f^*$, where $f^*,g^*$ are the nonincreasing 
rearrangements of $f,g$. For all~$R > 0$, we have
$$
   \int_{B_R} f^\#(x)\d x \,=\, \int_0^{c_N R^N} f^*(s)\d s~, \quad
   \int_{B_R} |x|^N f^\#(x)\d x \,=\, \frac{1}{c_N} \int_0^{c_N R^N} 
   s f^*(s)\d s~.
$$
Using the first equality, we can rewrite assumptions a) and c) as 
follows:\\[1mm]
\null\quad a') $\int_0^r h(s)\d s \ge 0$ for all $r > 0$, \\[1mm]
\null\quad c') $\int_0^\infty h(s)\d s = 0$.\\[1mm]
In view of \reff{concin} and the second equality, assumptions b) and
d) also imply\\[1mm]
\null\quad d') $\int_0^\infty s h(s)\d s \ge 0$, with equality if 
and only if $f = f^\#$. 

\medskip
Let $\cH(r) = \int_0^r h(s)\d s = - \int_r^\infty h(s)\d s \ge 0$. 
Clearly, $r\cH(r) \le \int_r^\infty s|h(s)|\d s \to 0$ as 
$r \to +\infty$, hence
$$
   \int_0^\infty s h(s)\d s + \int_0^\infty \cH(s)\d s 
   \,=\, \int_0^\infty (s h(s) + \cH(s))\d s  \,=\,
   r \cH(r)\Big|_{r=0}^{r=\infty} \,=\, 0~.
$$
By a') and d'), this implies that $\int_0^\infty s h(s)\d s = 0$ 
(hence $f = f^\#$) and $\cH \equiv 0$ (hence we have~$f^\# = g^\# = g$).
This concludes the proof. \QED

\medskip
Assume now that $f : \real^N \times [0,+\infty) \to \real$ is 
a solution of the convection-diffusion equation
\begin{equation}\label{CD}
  \partial_t f(x,t) + U(x,t) \cdot \nabla f(x,t) \,=\, 
  \Delta f(x,t)~, \quad x \in \real^N~, \quad t > 0~,
\end{equation}
where $U : \real^N \times [0,+\infty) \to \real^N$ is a smooth, 
divergence-free vector field which is bounded together with all 
its derivatives. Assume moreover that the initial condition $f_0(x) 
= f(x,0)$ is continuous and decays rapidly at infinity. 

\begin{proposition}\label{prop2}
Under the assumptions above, the solution of \reff{CD} with initial 
data $f_0$ satisfies
$$
   f(\cdot,t) \,\preceq\, \e^{t\Delta}f_0^\#~, \quad 
   \hbox{for all } t \ge 0 ~.
$$
\end{proposition}

\proof
The proof uses two main ingredients:\\[1mm]
i) If $f \preceq g$ and $g = g^\#$, then $\e^{t\Delta}f \preceq
\e^{t\Delta}g$ for all $t > 0$, see (\cite{ALT2}, Proposition~3).
This is because the integral kernel of the heat semigroup 
is positive, radially symmetric, and decreasing along rays (in 
other words, it coincides with its own rearrangement). \\[1mm]
ii) Let $S(t_1,t_2)$ be the evolution operator associated 
to the linear nonautonomous transport equation 
$$
  \partial_t f(x,t) + U(x,t) \cdot \nabla f(x,t) \,=\, 0~, 
  \quad x \in \real^N~, \quad t > 0~,
$$
so that $f(\cdot,t_1) = S(t_1,t_2) f(\cdot,t_2)$. Since $U(x,t)$ is 
globally Lipschitz  and divergence-free, one has~$S(t_1,t_2)f = f \circ 
\Psi(t_1,t_2)^{-1}$, where $\Psi(t_1,t_2) : \real^N \to \real^N$ is a 
measure preserving transformation. It follows immediately that 
$(S(t_1,t_2)f)^\# = f^\#$.  

\medskip
To prove the claim, we use an approximation procedure in the
spirit of Trotter's formula, see (\cite{taylor}, Section 11.A). 
Fix $t > 0$. Given $n \in \intplus^*$, we define finite sequences 
$\{f^n_k\}_{k=0,\dots,n}$ and~$\{g^n_k\}_{k=0,\dots,n}$ inductively by
$$
\begin{array}{rcll}
  f^n_{k+1} &=& {\displaystyle S\Bigl(\frac{(k{+}1)t}{n},
  \frac{kt}{n}\Bigr)\,\e^{\frac{t}{n}\Delta}f^n_k}~, \quad 
  & f^n_0 = f_0~,\\[2mm]
  g^n_{k+1} &=& {\displaystyle \e^{\frac{t}{n}\Delta}g^n_k}~, \quad 
  & g^n_0 = f_0^\#~.
\end{array}
$$
Using properties i) and ii) above, it is easy to verify by induction 
over $k$ that $f^n_k \preceq g^n_k$ for all~$k \in \{0,\dots,n\}$.
In particular,~$f^n_n \preceq g^n_n$. But~$g^n_n = \e^{t\Delta}
f_0^\#$, and a classical result  shows that~$f^n_n \to f(\cdot,t)$ in~$L^1(\real^N)$ as~$n \to \infty$,
so that~$(f^n_n)^\# \to f^\#(\cdot,t)$ also. Thus passing to the 
limit we obtain~$f(\cdot,t) \preceq \e^{t\Delta}f_0^\#$. \QED


\section{Uniqueness proof using symmetrization}
\label{unique2} 
\setcounter{equation}{0}

Let again $\omega(x,t)$ be a solution of \reff{Veq} satisfying 
the hypotheses of Theorem~\ref{main}, and let $u(x,t)$ be the 
associated velocity field. As usual we assume that $T = +\infty$ and 
that $\alpha = \inttwo \omega(x,t)\d x > 0$. For any $t > 0$, 
let $\omega^\#(x,t)$ be the symmetric nonincreasing rearrangement
of $\omega(x,t)$ with respect to the spatial variable $x \in \real^2$.
Using \reff{gaussbound} and the fact that the rearrangement is 
order preserving, we find
$$
  0 \,<\, \omega^\#(x,t) \,\le\, \frac{K_1 \alpha}{t}
  \,\e^{-\beta\frac{|x|^2}{4t}}~, \quad x \in \real^2~, \quad t > 0~.
$$
Since $\inttwo \omega^\#(x,t)\d x = \inttwo \omega(x,t)\d x = \alpha$,
we deduce that $\omega^\#(\cdot,t) \weakto \alpha\delta_0$ as 
$t \to 0+$. 

Fix $t > s > 0$. Applying Proposition~\ref{prop2} with 
$N = 2$, $f(x,t') = \omega(x,t'+s)$, and $U(x,t') = u(x,t'+s)$,
we obtain
$$
   \omega(\cdot,t) \,\preceq\, \e^{(t-s)\Delta}
   \omega^\#(\cdot,s)~.
$$
As $s \to 0+$, the right-hand side converges pointwise (hence
also in $L^1(\real^2)$ by the dominated convergence theorem) 
to the Oseen vortex $\Omega(\cdot,t)$ defined in \reff{Omegadef}.
It follows that
\begin{equation}\label{domtom}
   \omega(\cdot,t) \,\preceq\, \Omega(\cdot,t)~, \quad 
   \hbox{for all } t > 0~.
\end{equation}
We now fix $t > 0$ and apply Proposition~\ref{prop1} with $N = 2$, 
$f = \omega(\cdot,t)$, and $g = \Omega(\cdot,t)$. All assumptions
are satisfied due to \reff{mom0}, \reff{mom2}, and \reff{domtom}. 
We conclude that $\omega(x,t) = \Omega(x,t)$ for all $x \in \real^2$, 
which is the desired result. \QED


\bibliographystyle{plain}

\begin{thebibliography}{10}

\bibitem{ALT0}
A. Alvino, P.-L. Lions, and G. Trombetti. 
On optimization problems with prescribed rearrangements. 
{\em Nonlinear Anal. T.M.A.} {\bf 13} (1989), 185--220.

\bibitem{ALT1}
A. Alvino, P.-L. Lions, and G. Trombetti. 
Comparison results for elliptic and parabolic equations via 
Schwarz symmetrization. 
{\em Ann. Inst. H. Poincar\'e Anal. Non Lin\'eaire} {\bf 7} (1990), 37--65.

\bibitem{ALT2}
A. Alvino, P.-L. Lions, and G. Trombetti. 
Comparison results for elliptic and parabolic equations via 
symmetrization: a new approach. 
{\em Differential Integral Equations} {\bf 4} (1991), 25--50.

\bibitem{AMTU}
A. Arnold, P. Markowich, G. Toscani, and A. Unterreiter.
On convex {S}obolev inequalities and the rate of convergence to
equilibrium for Fokker-Planck type equations.
{\em Comm. Partial Differential Equations} {\bf 26} (2001), 43--100.

\bibitem{bandle}
C. Bandle.  
{\sl Isoperimetric inequalities and applications}. 
Monographs and Studies in Mathematics {\bf 7}. 
Pitman, London, 1980. 

\bibitem{ben-artzi}
M.~Ben-Artzi.
Global solutions of two-dimensional Navier-Stokes and Euler
equations.
{\em Arch. Rational Mech. Anal.} {\bf 128} (1994), 329--358.

\bibitem{brezis}
H.~Brezis.
Remarks on the preceding paper by M. Ben-{A}rtzi: ``Global
solutions of two-dimensional Navier-Stokes and Euler
equations''.
{\em Arch. Rational Mech. Anal.} {\bf 128} (1994), 359--360.

\bibitem{carlen-loss}
E.~A. Carlen and M.~Loss.
Optimal smoothing and decay estimates for viscously damped
conservation laws, with applications to the {$2$}-D Navier-Stokes
equation.
{\em Duke Math. J.} {\bf 81} (1995), 135--157 (1996). 

\bibitem{carpio}
A.~Carpio.
Asymptotic behavior for the vorticity equations in dimensions two and
three.
{\em Comm. Partial Differential Equations} {\bf 19} (1994), 827--872.

\bibitem{cottet}
G.-H. Cottet.
Equations de Navier-Stokes dans le plan avec tourbillon initial
mesure. 
{\em C. R. Acad. Sci. Paris S\'er. I Math.} {\bf 303} (1986), 105--108.

\bibitem{GG}
I. Gallagher and Th. Gallay.
Uniqueness for the two-dimensional Navier-Stokes equation with
a measure as initial vorticity. 
Preprint, Ecole Polytechnique, 2004. Available at 
{\tt http://www.arXiv.org/math.AP/0406297}.

\bibitem{gallay-wayne}
Th. Gallay and C.~E. Wayne.
Global stability of vortex solutions of the two-dimensional
Navier-Stokes equation. {\em Comm. Math. Phys.}, to appear.
Preprint version available at {\tt http://www.arXiv.org/math.AP/0402449}.

\bibitem{GMO}
Y.~Giga, T.~Miyakawa, and H.~Osada.
Two-dimensional Navier-Stokes flow with measures as initial
vorticity.
{\em Arch. Rational Mech. Anal.} {\bf 104} (1988), 223--250.

\bibitem{gigabook}
M.-H.~Giga and Y.~Giga.
{\sl Nonlinear partial differential equations -- asymptotic 
behaviour of solutions and self-similar solutions}.
Book in preparation. 

\bibitem{kato}
T.~Kato.
The Navier-Stokes equation for an incompressible fluid in
$\real^2$ with a measure as the initial vorticity.
{\em Differential Integral Equations} {\bf 7} (1994), 949--966.

\bibitem{lieb-loss}
E. Lieb and M. Loss. 
{\sl Analysis.} Graduate Studies in Mathematics {\bf 14}. 
American Mathematical Society, Providence, RI, 1997. 

\bibitem{osada}
H.~Osada.
Diffusion processes with generators of generalized divergence form.
{\em J. Math. Kyoto Univ.} {\bf 27} (1987), 597--619.

\bibitem{talenti}
G. Talenti.  
Elliptic equations and rearrangements. 
{\em Ann. Scuola Norm. Sup. Pisa} {\bf 3} (1976), 697--718.

\bibitem{taylor}
M.~Taylor.
{\sl Partial differential equations II. Qualitative studies
of linear equations}.
Applied Mathematical Sciences {\bf 116}, Springer, New-York, 1996. 

\end{thebibliography}

\end{document}